\documentclass{amsart}%
\usepackage{amsmath}
\usepackage{amssymb}
\usepackage{amsfonts}%
\usepackage{graphicx}
\providecommand{\U}[1]{\protect \rule{.1in}{.1in}}
\newtheorem{theorem}{Theorem}
\theoremstyle{plain}

\newtheorem{corollary}{Corollary}

\newtheorem{definition}{Definition}

\newtheorem{lemma}{Lemma}

\newtheorem{proposition}{Proposition}
\newtheorem{remark}{Remark}

\numberwithin{equation}{section}
\begin{document}
\title[vanishing generalized Morrey estimates]{vanishing generalized Morrey spaces and commutators of Marcinkiewicz integrals
with rough kernel associated with schr\"{o}dinger operator}
\author{F. GURBUZ}
\address{ANKARA UNIVERSITY, FACULTY OF SCIENCE, DEPARTMENT OF MATHEMATICS, TANDO\u{G}AN
06100, ANKARA, TURKEY }
\curraddr{}
\email{feritgurbuz84@hotmail.com}
\urladdr{}
\thanks{}
\thanks{}
\thanks{}
\date{}
\dedicatory{ }
\begin{abstract}
Let $L=-\Delta+V\left(  x\right)  $ be a Schr\"{o}dinger operator, where
$\Delta$ is the Laplacian on ${\mathbb{R}^{n}}$, while nonnegative potential
$V\left(  x\right)  $ belonging to the reverse H\"{o}lder class. We establish
the boundedness of the commutators of Marcinkiewicz integrals with rough
kernel associated with schr\"{o}dinger operator on vanishing generalized
Morrey spaces.

\end{abstract}
\maketitle

\section{Introduction and main results}

Because of the need for study of the local behavior of solutions of second
order elliptic partial differential equations (PDEs) and together with the now
well-studied Sobolev Spaces, constitude a formidable three parameter family of
spaces useful for proving regularity results for solutions to various PDEs,
especially for non-linear elliptic systems, in 1938, Morrey \cite{Morrey}
introduced the classical Morrey spaces $L_{p,\lambda}$ which naturally are
generalizations of the classical Lebesgue spaces.

We will say that a function $f\in L_{p,\lambda}=L_{p,\lambda}\left(
{\mathbb{R}^{n}}\right)  $ if%
\begin{equation}
\sup_{x\in{\mathbb{R}^{n},r>0}}\left[  r^{-\lambda}%
{\displaystyle \int \limits_{B(x,r)}}
\left \vert f\left(  y\right)  \right \vert ^{p}dy\right]  ^{1/p}<\infty
.\label{1.1}%
\end{equation}
Here, $1<p<\infty$ and $0<\lambda<n$ and the quantity of (\ref{1.1}) is the
$\left(  p,\lambda \right)  $-Morrey norm, denoted by $\left \Vert f\right \Vert
_{L_{p,\lambda}}$. In recent years, more and more researches focus on function
spaces based on Morrey spaces to fill in some gaps in the theory of Morrey
type spaces (see, for example, \cite{Gurbuz1, Gurbuz2, Gurbuz3, Gurbuz4,
Gurbuz5, Pal}). Moreover, these spaces are proved useful in harmonic analysis
and PDEs. But, this topic exceeds the scope of this paper. Thus, we omit the
details here. On the other hand, the study of the operators of harmonic
analysis in vanishing Morrey space, in fact has been almost not touched. A
version of the classical Morrey space $L_{p,\lambda}({\mathbb{R}^{n}})$ where
it is possible to approximate by "nice" functions is the so called vanishing
Morrey space $VL_{p,\lambda}({\mathbb{R}^{n}})$ has been introduced by Vitanza
in \cite{Vitanza1} and has been applied there to obtain a regularity result
for elliptic PDEs. This is a subspace of functions in $L_{p,\lambda
}({\mathbb{R}^{n}})$, which satisfies the condition%
\[
\lim_{r\rightarrow0}\sup_{\underset{0<t<r}{x\in{\mathbb{R}^{n}}}}\left[
t^{-\lambda}%
{\displaystyle \int \limits_{B(x,t)}}
\left \vert f\left(  y\right)  \right \vert ^{p}dy\right]  ^{1/p}=0,
\]
where $1<p<\infty$ and $0<\lambda<n$ for brevity, so that%
\[
VL_{p,\lambda}({\mathbb{R}^{n}})=\left \{  f\in L_{p,\lambda}({\mathbb{R}^{n}%
}):\lim_{r\rightarrow0}\sup_{\underset{0<t<r}{x\in{\mathbb{R}^{n}}}}%
t^{-\frac{\lambda}{p}}\Vert f\Vert_{L_{p}(B(x,t))}=0\right \}  .
\]
Later in \cite{Vitanza2} Vitanza has proved an existence theorem for a
Dirichlet problem, under weaker assumptions than in \cite{Miranda} and a
$W^{3,2}$ regularity result assuming that the partial derivatives of the
coefficients of the highest and lower order terms belong to vanishing Morrey
spaces depending on the dimension. For the properties and applications of
vanishing Morrey spaces, see also \cite{Cao-Chen}.

After studying Morrey spaces in detail, researchers have passed to the concept
of generalized Morrey spaces. Firstly, motivated by the work of \cite{Morrey},
Mizuhara \cite{Miz} introduced generalized Morrey spaces $M_{p,\varphi}$ as follows:

\begin{definition}
\cite{Miz} \textbf{(generalized Morrey space)\ }Let $\varphi(x,r)$ be a
positive measurable function on ${\mathbb{R}^{n}}\times(0,\infty)$. If
$0<p<\infty$, then the generalized Morrey space $M_{p,\varphi}\equiv
M_{p,\varphi}({\mathbb{R}^{n}})$ is defined by%
\[
\left \{  f\in L_{p}^{loc}({\mathbb{R}^{n}}):\Vert f\Vert_{M_{p,\varphi}}%
=\sup \limits_{x\in{\mathbb{R}^{n}},r>0}\varphi(x,r)^{-1}\Vert f\Vert
_{L_{p}(B(x,r))}<\infty \right \}  .
\]

\end{definition}

Obviously, the above definition recover the definition of $L_{p,\lambda
}({\mathbb{R}^{n}})$ if we choose $\varphi(x,r)=r^{\frac{\lambda}{p}}$, that
is
\[
L_{p,\lambda}\left(  {\mathbb{R}^{n}}\right)  =M_{p,\varphi}\left(
{\mathbb{R}^{n}}\right)  \mid_{\varphi(x,r)=r^{\frac{\lambda}{p}}}.
\]

Everywhere in the sequel we assume that $\inf \limits_{x\in{\mathbb{R}^{n}%
},r>0}\varphi(x,r)>0$ which makes the above spaces non-trivial, since the
spaces of bounded functions are contained in these spaces. We point out that
$\varphi(x,r)$ is a measurable non-negative function and no monotonicity type
condition is imposed on these spaces.

Throughout the paper we assume that $x\in{\mathbb{R}^{n}}$ and $r>0$ and also
let $B(x,r)$ denotes the open ball centered at $x$ of radius $r$, $B^{C}(x,r)$
denotes its complement and $|B(x,r)|$ is the Lebesgue measure of the ball
$B(x,r)$ and $|B(x,r)|=v_{n}r^{n}$, where $v_{n}=|B(0,1)|$.

Now, recall that the concept of the vanishing generalized Morrey spaces
$VM_{p,\varphi}({\mathbb{R}^{n}})$ has been introduced in \cite{N. Samko}.

\begin{definition}
\cite{N. Samko} \textbf{(vanishing generalized Morrey space) }Let
$\varphi(x,r)$ be a positive measurable function on ${\mathbb{R}^{n}}%
\times(0,\infty)$ and $1\leq p<\infty$. The vanishing generalized Morrey space
$VM_{p,\varphi}({\mathbb{R}^{n}})$ is defined as the spaces of functions $f\in
L_{p}^{loc}({\mathbb{R}^{n}})$ such that%
\begin{equation}
\lim \limits_{r\rightarrow0}\sup \limits_{x\in{\mathbb{R}^{n}}}\varphi
(x,r)^{-1}\, \int_{B(x,r)}|f(y)|^{p}dy=0.\label{1*}%
\end{equation}

\end{definition}

Everywhere in the sequel we assume that%
\begin{equation}
\lim_{t\rightarrow0}\frac{t^{\frac{n}{p}}}{\varphi(x,t)}=0,\label{2}%
\end{equation}
and%
\begin{equation}
\sup_{0<t<\infty}\frac{t^{\frac{n}{p}}}{\varphi(x,t)}<\infty,\label{3}%
\end{equation}
which make the spaces $VM_{p,\varphi}({\mathbb{R}^{n}})$ non-trivial, because
bounded functions with compact support belong to this space. The spaces
$VM_{p,\varphi}({\mathbb{R}^{n}})$ and $WVM_{p,\varphi}({\mathbb{R}^{n}})$ are
Banach spaces with respect to the norm (see, for example \cite{N. Samko})%

\begin{equation}
\Vert f\Vert_{VM_{p,\varphi}}=\sup \limits_{x\in{\mathbb{R}^{n}},r>0}%
\varphi(x,r)^{-1}\Vert f\Vert_{L_{p}(B(x,r))},\label{4}%
\end{equation}%
\begin{equation}
\Vert f\Vert_{WVM_{p,\varphi}}=\sup \limits_{x\in{\mathbb{R}^{n}},r>0}%
\varphi(x,r)^{-1}\Vert f\Vert_{WL_{p}(B(x,r)),}\label{5}%
\end{equation}
respectively. For the properties and applications of vanishing generalized
Morrey spaces, see also \cite{Akbulut-Kuzu}. In \cite{Akbulut-Kuzu}, the
boundedness of the Marcinkiewicz integrals with rough kernel associated with
schr\"{o}dinger operator on vanishing generalized Morrey spaces $VM_{p,\varphi
}({\mathbb{R}^{n}})$ has been investigated.

On the other hand, suppose that $S^{n-1}$ is the unit sphere in ${\mathbb{R}%
^{n}}$ $(n\geq2)$ equipped with the normalized Lebesgue measure $d\sigma
=d\sigma \left(  x^{\prime}\right)  $.

In \cite{Stein58}, Stein has defined the Marcinkiewicz integral for higher
dimensions. Suppose that $\Omega$ satisfies the following conditions.

(a) $\Omega$ is the homogeneous function of degree zero on ${\mathbb{R}^{n}%
}\setminus \{0\}$, that is,
\begin{equation}
\Omega(\mu x)=\Omega(x),~~\text{for any}~~\mu>0,x\in{\mathbb{R}^{n}}%
\setminus \{0\}.\label{1}%
\end{equation}

(b) $\Omega$ has mean zero on $S^{n-1}$, that is,
\begin{equation}
\int \limits_{S^{n-1}}\Omega(x^{\prime})d\sigma(x^{\prime})=0,\label{2*}%
\end{equation}
where $x^{\prime}=\frac{x}{\left \vert x\right \vert }$ for any $x\neq0$.

(c) $\Omega \in Lip_{\gamma}(S^{n-1})$, $0<\gamma \leq1$, that is there exists a
constant $M>0$ such that,
\[
|\Omega(x^{\prime})-\Omega(y^{\prime})|\leq M|x^{\prime}-y^{\prime}|^{\gamma
}~~\text{for any}~~x^{\prime},y^{\prime}\in S^{n-1}.
\]

(d) $\Omega \in L_{1}(S^{n-1})$.

The Marcinkiewicz integral operator of higher dimension $\mu_{\Omega}$ is
defined by
\[
\mu_{\Omega}(f)(x)=\left(  \int \limits_{0}^{\infty}|F_{\Omega,t}%
(f)(x)|^{2}\frac{dt}{t^{3}}\right)  ^{1/2},
\]
where
\[
F_{\Omega,t}(f)(x)=\int \limits_{|x-y|\leq t}\frac{\Omega(x-y)}{|x-y|^{n-1}%
}f(y)dy.
\]

Since Stein's work in 1958, the continuity of Marcinkiewicz integral has been
extensively studied as a research topic and also provides useful tools in
harmonic analysis \cite{St, Stein93, Torch}.

\begin{remark}
We easily see that the Marcinkiewicz integral operator of higher dimension
$\mu_{\Omega}$ can be regarded as a generalized version of the classical
Marcinkiewicz integral in the one dimension case. Also, it is easy to see that
$\mu_{\Omega}$ is a special case of the Littlewood-Paley $g$-function if we
take%
\[
g\left(  x\right)  =\Omega \left(  x^{\prime}\right)  \left \vert x\right \vert
^{-n+1}\chi_{\left \vert x\right \vert \leq1}\left(  \left \vert x\right \vert
\right)  .
\]

\end{remark}

When $\Omega$ satisfies some size conditions, the kernel of the operator
$\mu_{\Omega}$ has no regularity, and so the operator $\mu_{\Omega}$ is called
rough Marcinkiewicz integral operator.{\ }The theory of Operators with rough
kernel is a well studied area (see \cite{Gurbuz1, Gurbuz3, Gurbuz4, Gurbuz5}
for example).

For simplicity of notation, $\Omega$ is always homogeneous function of degree
zero and satisfies
\[
\Omega \in L_{q}(S^{n-1}),\qquad1<q\leq \infty
\]
and (\ref{2*}) throughout this paper if there are no special instructions.

Now we give the definition of the commutator generalized by $\mu_{\Omega}$ and
$b$ by%
\[
\mu_{\Omega,b}(f)(x)=\left(
{\displaystyle \int \limits_{0}^{\infty}}
|F_{\Omega,t,b}(f)(x)|^{2}\frac{dt}{t^{3}}\right)  ^{1/2},
\]
where%
\[
F_{\Omega,t,b}(f)(x)=%
{\displaystyle \int \limits_{|x-y|\leq t}}
\frac{\Omega(x-y)}{|x-y|^{n-1}}[b(x)-b(y)]f(y)dy.
\]

Let $f\in L_{1}^{loc}({\mathbb{R}^{n}})$. The rough Hardy-Littlewood maximal
operator $M_{\Omega}$ and commutator of the Hardy-Littlewood maximal operator
with rough kernel are defined by%

\[
M_{\Omega}f\left(  x\right)  =\sup_{t>0}\frac{1}{\left \vert B\left(
x,t\right)  \right \vert }\int \limits_{B\left(  x,t\right)  }\left \vert
\Omega \left(  x-y\right)  \right \vert \left \vert f\left(  y\right)
\right \vert dy,
\]%
\[
M_{\Omega,b}\left(  f\right)  (x)=\sup_{t>0}|B(x,t)|^{-1}\int \limits_{B(x,t)}%
\left \vert b\left(  x\right)  -b\left(  y\right)  \right \vert \left \vert
\Omega \left(  x-y\right)  \right \vert |f(y)|dy,
\]
respectively.

The following results concerning the boundedness of commutator operators
$\mu_{\Omega,b}$ and $M_{\Omega,b}$ on $L_{p}$ space are known.

\begin{theorem}
\label{teo7}(see \cite{Ding-Lu-Yabuta}) Let $1\leq p<\infty$, $\Omega \in
L_{q}(S^{n-1})$, $1<q\leq \infty$ satisfies (\ref{1}), (\ref{2*}) and $b\in
BMO\left(  {\mathbb{R}^{n}}\right)  $. Then, for $p>1$ $\mu_{\Omega,b}$ is
bounded on $L_{p}\left(  {\mathbb{R}^{n}}\right)  $ and for $p=1$ from
$L_{1}\left(  {\mathbb{R}^{n}}\right)  $ to $WL_{1}\left(  {\mathbb{R}^{n}%
}\right)  $.
\end{theorem}

\begin{theorem}
\label{teo8}(see \cite{Al-Bag-Kurt-Perez}) Let $1<p<\infty$, $\Omega \in
L_{q}(S^{n-1})$, $1<q\leq \infty$ satisfies (\ref{1}) and $b\in BMO\left(
{\mathbb{R}^{n}}\right)  $. Then, for every $q^{\prime}<p<\infty$ or $1<p<q$,
there is a constant $C$ independent of $f$ such that
\[
\left \Vert M_{\Omega,b}\left(  f\right)  \right \Vert _{L_{p}}\leq C\left \Vert
f\right \Vert _{L_{p}}.
\]
Moreover, for $p>1$ $M_{\Omega,b}$ is bounded on $L_{p}\left(  {\mathbb{R}%
^{n}}\right)  $ and for $p=1$ from $L_{1}\left(  {\mathbb{R}^{n}}\right)  $ to
$WL_{1}\left(  {\mathbb{R}^{n}}\right)  $.
\end{theorem}

On the other hand, in this paper we consider the Schr\"{o}dinger operator%
\[
L=-\Delta+V\left(  x\right)  \text{ on }{\mathbb{R}^{n}},\qquad n\geq3
\]
where $V\left(  x\right)  $ is a nonnegative potential belonging to the
reverse H\"{o}lder class $RH_{q}$, for some exponent $q\geq \frac{n}{2}$; that
is, there exists a constant $C$ such that the reverse H\"{o}lder inequality%
\begin{equation}
\left(  \frac{1}{\left \vert B\right \vert }%
{\displaystyle \int \limits_{B}}
V\left(  x\right)  ^{q}dx\right)  ^{\frac{1}{q}}\leq \frac{C}{\left \vert
B\right \vert }%
{\displaystyle \int \limits_{B}}
V\left(  x\right)  dx,\label{3-}%
\end{equation}
holds for every ball $B\subset{\mathbb{R}^{n}}$; see \cite{Shen}.

We introduce the definition of the reverse H\"{o}lder index of $V$ as
$q_{0}=\sup \left \{  q:V\in RH_{q}\right \}  $. It is worth pointing out that
the $RH_{q}$ class is that, if $V\in RH_{q}$ for some $q>1$, then there exists
$\varepsilon>0$, which depends only on $n$ and the constant $C$ in (\ref{3-}),
such that $V\in$ $RH_{q+\varepsilon}$. Therefore, under the assumption $V\in
RH_{\frac{n}{2}}$, we may conclude $q_{0}>\frac{n}{2}$. Throughout this paper,
we always assume that $0\neq V\in RH_{n}$. In particular, Shen \cite{Shen} has
considered $L_{p}$ estimates for Schr\"{o}dinger operators $L$ with certain
potentials which include Schr\"{o}dinger Riesz transforms $R_{j}^{L}%
=\frac{\partial}{\partial x_{j}}L^{-\frac{1}{2}}$, $j=1,\cdots,n$. Then,
Dziuba\'{n}ski and Zienkiewicz \cite{Dzi-Zien} has introduced the Hardy type
space $H_{L}^{1}\left(  {\mathbb{R}^{n}}\right)  $ associated with the
Schr\"{o}dinger operator $L$, which is larger than the classical Hardy space
$H^{1}\left(  {\mathbb{R}^{n}}\right)  $.

Similar to the Marcinkiewicz integral operator with rough kernel $\mu_{\Omega
}$, we define the Marcinkiewicz integral operator with rough kernel
$\mu_{j,\Omega}^{L}$ associated with the Schr\"{o}dinger operator $L$ by%
\[
\mu_{j,\Omega}^{L}f\left(  x\right)  =\left(
{\displaystyle \int \limits_{0}^{\infty}}
\left \vert
{\displaystyle \int \limits_{\left \vert x-y\right \vert \leq t}}
\left \vert \Omega \left(  x-y\right)  \right \vert K_{j}^{L}\left(  x,y\right)
f\left(  y\right)  dy\right \vert ^{2}\frac{dt}{t^{3}}\right)  ^{\frac{1}{2}},
\]
where $K_{j}^{L}\left(  x,y\right)  =\widetilde{K_{j}^{L}}\left(  x,y\right)
\left \vert x-y\right \vert $ and $\widetilde{K_{j}^{L}}\left(  x,y\right)  $ is
the kernel of $R_{j}=\frac{\partial}{\partial_{x_{j}}}L^{-\frac{1}{2}}$,
$j=1,\ldots,n$. In particular, when $V=0$, $K_{j}^{\Delta}\left(  x,y\right)
=\widetilde{K_{j}^{\Delta}}\left(  x,y\right)  \left \vert x-y\right \vert
=\left(  \left(  x_{j}-y_{j}\right)  /\left \vert x-y\right \vert \right)
/\left \vert x-y\right \vert ^{n-1}$ and $\widetilde{K_{j}^{\Delta}}\left(
x,y\right)  $ is the kernel of $R_{j}=\frac{\partial}{\partial_{x_{j}}}%
\Delta^{-\frac{1}{2}}$, $j=1,\ldots,n$. In this paper, we write $K_{j}\left(
x,y\right)  =K_{j}^{\Delta}\left(  x,y\right)  $ and $\mu_{j,\Omega}%
=\mu_{j,\Omega}^{\Delta}$ and so, $\mu_{j,\Omega}^{\Delta}$ is defined by%
\[
\mu_{j,\Omega}f\left(  x\right)  =\left(
{\displaystyle \int \limits_{0}^{\infty}}
\left \vert
{\displaystyle \int \limits_{\left \vert x-y\right \vert \leq t}}
\left \vert \Omega \left(  x-y\right)  \right \vert K_{j}\left(  x,y\right)
f\left(  y\right)  dy\right \vert ^{2}\frac{dt}{t^{3}}\right)  ^{\frac{1}{2}}.
\]

Obviously, $\mu_{j}$ are classical Marcinkiewicz functions. Therefore, it will
be an interesting thing to study the properties of $\mu_{j,\Omega}^{L} $.

Given an operator $\mu_{j,\Omega}^{L}$, and a function $b$, we define the
commutator of $\mu_{j,\Omega}^{L}$ and $b$ by%
\[
\mu_{j,\Omega,b}^{L}f\left(  x\right)  =[b,\mu_{j,\Omega}^{L}]f(x)=b(x)\,
\mu_{j,\Omega}^{L}f(x)-\mu_{j,\Omega}^{L}(bf)(x).
\]
If $\mu_{j,\Omega}^{L}$ is defined by integration against a kernel for certain
$x$, such as when $\mu_{j,\Omega}^{L}$ is Marcinkiewicz integral operator with
rough kernel associated with the Schr\"{o}dinger operator $L$, we have that
this becomes%
\[
\mu_{j,\Omega,b}^{L}f\left(  x\right)  =[b,\mu_{j,\Omega}^{L}]f(x)=\left(
{\displaystyle \int \limits_{0}^{\infty}}
\left \vert
{\displaystyle \int \limits_{\left \vert x-y\right \vert \leq t}}
\left \vert \Omega \left(  x-y\right)  \right \vert K_{j}^{L}\left(  x,y\right)
\left[  b\left(  x\right)  -b\left(  y\right)  \right]  f\left(  y\right)
dy\right \vert ^{2}\frac{dt}{t^{3}}\right)  ^{\frac{1}{2}},
\]
for all $x$ for which the integral representation of $\mu_{j,\Omega}^{L}$
holds. It is worth noting that for a constant $C$, if $\mu_{j,\Omega}^{L}$ is
linear we have,%
\begin{align*}
\lbrack b+C,\mu_{j,\Omega}^{L}]f  & =\left(  b+C\right)  \mu_{j,\Omega}%
^{L}f-\mu_{j,\Omega}^{L}(\left(  b+C\right)  f)\\
& =b\mu_{j,\Omega}^{L}f+C\mu_{j,\Omega}^{L}f-\mu_{j,\Omega}^{L}\left(
bf\right)  -C\mu_{j,\Omega}^{L}f\\
& =[b,\mu_{j,\Omega}^{L}]f.
\end{align*}
This leads one to intuitively look to spaces for which we identify functions
which differ by constants, and so it is no surprise that $b\in BMO$ (bounded
mean oscillation {space) has had the most historical significance.}

Now, for a given potential $V\in RH_{q}$, with $q>\frac{n}{2}$, we introduce
the auxiliary function
\[
\rho \left(  x\right)  =\frac{1}{m_{_{V}}\left(  x\right)  }=\sup_{r>0}\left \{
r:\frac{1}{r^{n-2}}%
{\displaystyle \int \limits_{B\left(  x,r\right)  }}
V\left(  y\right)  dy\leq1\right \}  ,\qquad x\in{\mathbb{R}^{n}.}%
\]

The above assumptions $\rho \left(  x\right)  $ are finite, for all
$x\in{\mathbb{R}^{n}}$. Obviously, $0<m_{_{V}}\left(  x\right)  <\infty$ if
$V\neq0 $. In particular, $m_{_{V}}\left(  x\right)  =1$ with $V=1$ and
$m_{_{V}}\left(  x\right)  \sim \left(  1+\left \vert x\right \vert \right)  $
with $V=\left \vert x\right \vert ^{2}$.

\begin{proposition}
(see \cite{Shen}) There exist $C$ and $k_{0}\geq1$ such that
\[
C^{-1}\rho \left(  x\right)  \left(  1+\frac{\left \vert x-y\right \vert }%
{\rho \left(  x\right)  }\right)  ^{-k_{0}}\leq \rho \left(  y\right)  \leq
C\rho \left(  x\right)  \left(  1+\frac{\left \vert x-y\right \vert }{\rho \left(
x\right)  }\right)  ^{\frac{k_{0}}{1+k_{0}}},
\]
for all $x$, $y\in{\mathbb{R}^{n}}$.

In particular, $\rho \left(  x\right)  \sim \rho \left(  y\right)  $, if
$\left \vert x-y\right \vert <C\rho \left(  x\right)  $. A ball $B\left(
x,\rho \left(  x\right)  \right)  $ is called critical.
\end{proposition}

\begin{proposition}
(see \cite{Dzi-Zien}) There exist a sequence of points $\left \{
x_{k}\right \}  _{k=1}^{\infty}$ in ${\mathbb{R}^{n}}$, so that the family
$B_{k}=B\left(  x_{k},\rho \left(  x_{k}\right)  \right)  $, $k\geq1$,
satisfies the following:

(1) $%
{\displaystyle \bigcup \limits_{k}}
B_{k}={\mathbb{R}^{n};}$

(2) There exist $N$ such that, for every $k\in$ $N$, card $\left \{
j:4B_{j}\cap4B_{k}\neq \emptyset \right \}  \leq N$.

\begin{lemma}
(see \cite{Tang-Dong}) For any $l>0$, there exists $C_{l}>0$ such that
\[
K_{j}^{L}\left(  x,y\right)  \leq \frac{C_{l}}{\left(  1+\frac{\left \vert
x-y\right \vert }{\rho \left(  y\right)  }\right)  ^{l}}\frac{1}{\left \vert
x-y\right \vert ^{n-1}},
\]
and%
\[
\left \vert K_{j}^{L}\left(  x,y\right)  -K_{j}\left(  x,y\right)  \right \vert
\leq C\frac{\rho \left(  x\right)  }{\left \vert x-y\right \vert ^{n-2}},
\]
where $\rho$ is the auxiliary function.
\end{lemma}
\end{proposition}

Tang and Dong \cite{Tang-Dong} have shown that Marcinkiewicz integral $\mu
_{j}^{L}$ is bounded on $L_{p}({\mathbb{R}^{n}})$, for $1<p<\infty$, and are
bounded from $L_{1}({\mathbb{R}^{n}})$ to $WL_{1}({\mathbb{R}^{n}})$.

Shen \cite{Shen} has given the following kernel estimate that we need.

\begin{lemma}
\label{lemma1}If $V\in RH_{n}$, then, one has

(i) for every $N$ there exist a constant $C$ such that%
\[
\left \vert K_{j}^{L}\left(  x,z\right)  \right \vert \leq \frac{C\left(
1+\frac{\left \vert x-z\right \vert }{\rho \left(  x\right)  }\right)  ^{-N}%
}{\left \vert x-z\right \vert ^{n-1}},
\]

(ii) for every $N$ and $0<\delta<\min \left \{  1,1-\frac{n}{q_{0}}\right \}  $,
there exists a constant $C$ such that%
\[
\left \vert K_{j}^{L}\left(  x,z\right)  -K_{j}^{L}\left(  y,z\right)
\right \vert \leq \frac{C\left \vert x-y\right \vert ^{\delta}\left(
1+\frac{\left \vert x-z\right \vert }{\rho \left(  x\right)  }\right)  ^{-N}%
}{\left \vert x-z\right \vert ^{n-1+\delta}},
\]
where $\left \vert x-y\right \vert <\frac{2}{3}\left \vert x-z\right \vert $,

(iii) if $K$ denotes the ${\mathbb{R}^{n}}$ vector valued kernel of the
classical Riesz operator, for every $0<\delta<2-\frac{n}{q_{0}}$, we have%
\[
\left \vert K_{j}^{L}\left(  x,z\right)  -K_{j}\left(  x,z\right)  \right \vert
\leq \frac{C}{\left \vert x-z\right \vert ^{n-1}}\left(  \frac{\left \vert
x-z\right \vert }{\rho \left(  z\right)  }\right)  ^{\delta},
\]
where $K_{j}\left(  x,z\right)  =K\left(  x,z\right)  \left \vert
x-z\right \vert $.
\end{lemma}

Inspired by \cite{Akbulut-Kuzu}, we give $BMO$ estimates for commutators of
Marcinkiewicz integrals with rough kernel associated with schr\"{o}dinger
operator on vanishing generalized Morrey spaces $VM_{p,\varphi}({\mathbb{R}%
^{n}})$.

We now make some conventions. Throughout this paper, we use the symbol
$A\lesssim B$ to denote that there exists a positive consant $C$ such that
$A\leq CB$. If $A\lesssim B$ and $B\lesssim A$, we then write $A\approx B$ and
say that $A$ and $B$ are equivalent. For a fixed $p\in \left[  1,\infty \right)
$, $p^{\prime}$ denotes the dual or conjugate exponent of $p$, namely,
$p^{\prime}=\frac{p}{p-1}$ and we use the convention $1^{\prime}=\infty$ and
$\infty^{\prime}=1$.

Our main results can be formulated as follows.

\begin{theorem}
\label{teo9*}Let $1<p<\infty$, $\Omega \in L_{q}(S^{n-1})$, $1<q\leq \infty$
satisfies (\ref{1}). Also, let $V\in RH_{n}$ and $b\in BMO\left(
{\mathbb{R}^{n}}\right)  $. Then, for every $q^{\prime}<p<\infty$ or $1<p<q$,
there is a constant $C$ independent of $f$ such that
\[
\left \Vert \mu_{j,\Omega,b}^{L}f\right \Vert _{L_{p}}\leq C\left \Vert
f\right \Vert _{L_{p}}.
\]

\end{theorem}

\begin{theorem}
\label{teo10*}Let $x_{0}\in{\mathbb{R}^{n}}$,$1<p<\infty$ and $b\in BMO\left(
{\mathbb{R}^{n}}\right)  $. Let $\Omega \in L_{q}(S^{n-1})$, $1<q\leq \infty$
satisfies (\ref{1}) and $V\in RH_{n}$. Then, for $q^{\prime}\leq p$ the
inequality
\begin{equation}
\Vert \mu_{j,\Omega,b}^{L}f\Vert_{L_{p}(B(x_{0},r))}\lesssim \Vert b\Vert_{\ast
}\,r^{\frac{n}{p}}\int \limits_{2r}^{\infty}\left(  1+\ln \frac{t}{r}\right)
t^{-\frac{n}{p}-1}\Vert f\Vert_{L_{p}(B(x_{0},t))}dt\label{47}%
\end{equation}
holds for any ball $B(x_{0},r)$ and for all $f\in L_{p}^{loc}({\mathbb{R}^{n}%
})$.

Also, for $p<q$ the inequality%
\begin{equation}
\Vert \mu_{j,\Omega,b}^{L}f\Vert_{L_{p}(B(x_{0},r))}\lesssim \Vert b\Vert_{\ast
}\,r^{\frac{n}{p}-\frac{n}{q}}\int \limits_{2r}^{\infty}\left(  1+\ln \frac
{t}{r}\right)  t^{\frac{n}{q}-\frac{n}{p}-1}\Vert f\Vert_{L_{p}(B(x_{0}%
,t))}dt\label{48}%
\end{equation}
holds for any ball $B(x_{0},r)$ and for all $f\in L_{p}^{loc}({\mathbb{R}^{n}%
})$.
\end{theorem}

\begin{theorem}
\label{teo11}Let $\Omega \in L_{q}(S^{n-1})$,$1<q\leq \infty$, satisfies
(\ref{1}) and $V\in RH_{n}$ . Let $1<p<\infty$ and $b\in BMO\left(
\mathbb{R}
^{n}\right)  $. For $q^{\prime}\leq p$ if the pair $(\varphi_{1},\varphi_{2})$
satisfies conditions (\ref{2})-(\ref{3}) and
\begin{equation}
c_{\delta}:=%
{\displaystyle \int \limits_{\delta}^{\infty}}
\left(  1+\ln \frac{t}{r}\right)  \sup_{x\in{\mathbb{R}^{n}}}\varphi_{1}\left(
x,t\right)  t^{-\frac{n}{p}-1}dt<\infty \label{6*}%
\end{equation}
for every $\delta>0$, and
\begin{equation}
\int \limits_{r}^{\infty}\left(  1+\ln \frac{t}{r}\right)  \frac{\varphi
_{1}(x,t)}{t^{\frac{n}{p}+1}}dt\leq C_{0}\frac{\varphi_{2}(x,r)}{r^{\frac
{n}{p}}},\label{7*}%
\end{equation}
and for $p<q$ if the pair $(\varphi_{1},\varphi_{2})$ satisfies conditions
(\ref{2})-(\ref{3}) and also%
\begin{equation}
c_{\delta^{\prime}}:=%
{\displaystyle \int \limits_{\delta^{\prime}}^{\infty}}
\left(  1+\ln \frac{t}{r}\right)  \sup_{x\in{\mathbb{R}^{n}}}\varphi_{1}\left(
x,t\right)  t^{-\frac{n}{p}+\frac{n}{q}-1}dt<\infty \label{8*}%
\end{equation}
for every $\delta^{\prime}>0$, and%
\begin{equation}
\int \limits_{r}^{\infty}\left(  1+\ln \frac{t}{r}\right)  \frac{\varphi
_{1}(x,t)}{t^{\frac{n}{p}-\frac{n}{q}+1}}dt\leq C_{0}\, \frac{\varphi
_{2}(x,r)}{r^{\frac{n}{p}-\frac{n}{q}}},\label{9*}%
\end{equation}
where $C_{0}$ does not depend on $x\in{\mathbb{R}^{n}}$ and $r>0$, then the
operators $\mu_{j,\Omega,b}^{L}$, $j=1,\ldots,n$ are bounded from
$VM_{p,\varphi_{1}}$ to $VM_{p,\varphi_{2}}$. Moreover,%
\begin{equation}
\left \Vert \mu_{j,\Omega,b}^{L}f\right \Vert _{VM_{p,\varphi_{2}}}%
\lesssim \left \Vert b\right \Vert _{\ast}\left \Vert f\right \Vert _{VM_{p,\varphi
_{1}}}.\label{10*}%
\end{equation}

\end{theorem}

\section{Some preliminaries}

We begin with some properties of $BMO\left(  {\mathbb{R}^{n}}\right)  $ spaces
which play a great role in the proofs of our main results.

Let us recall the defination of the space of $BMO({\mathbb{R}^{n}})$.

\begin{definition}
\label{definition}\cite{John-Nirenberg} The space $BMO({\mathbb{R}^{n}})$ of
functions of bounded mean oscillation consists of locally summable functions
with finite semi-norm%
\begin{equation}
\Vert b\Vert_{\ast}\equiv \Vert b\Vert_{BMO}=\sup_{x\in{\mathbb{R}^{n}}%
,r>0}\frac{1}{|B(x,r)|}%
{\displaystyle \int \limits_{B(x,r)}}
|b(y)-b_{B(x,r)}|dy<\infty,\label{0*}%
\end{equation}
where $b_{B(x,r)}$ is the mean value of the function $b$ on the ball $B(x,r)$.
The fact that precisely the mean value $b_{B(x,r)}$ figures in (\ref{0*}) is
inessential and one gets an equivalent seminorm if $b_{B(x,r)}$ is replaced by
an arbitrary constant $c:$%
\begin{equation}
\Vert b\Vert_{\ast}\approx \sup_{r>0}\inf_{c\in%
\mathbb{C}
}\frac{1}{\left \vert B\left(  x,r\right)  \right \vert }\int \limits_{B\left(
x,r\right)  }\left \vert b\left(  y\right)  -c\right \vert dy.\label{0**}%
\end{equation}
Indeed, it is obvious that (\ref{0*}) implies (\ref{0**}). If (\ref{0**})
holds, then
\[
\left \vert b_{B(x,r)}-c\right \vert =\left \vert \frac{1}{|B(x,r)|}%
\int \limits_{B\left(  x,r\right)  }\left(  b\left(  y\right)  -c\right)
dy\right \vert \leq C,
\]
so%
\[
\frac{1}{|B(x,r)|}%
{\displaystyle \int \limits_{B(x,r)}}
|b(y)-b_{B(x,r)}|dy\leq \frac{1}{|B(x,r)|}%
{\displaystyle \int \limits_{B(x,r)}}
\left(  \left \vert b\left(  y\right)  -c\right \vert +\left \vert c-b_{B(x,r)}%
\right \vert \right)  dy\leq2C.
\]

\end{definition}

Each bounded function $b\in BMO$. Moreover, $BMO$ contains unbounded
functions, in fact log$\left \vert x\right \vert $ belongs to $BMO$ but is not
bounded, so $L_{\infty}({\mathbb{R}^{n}})\subset BMO({\mathbb{R}^{n}})$.

In 1961 John and Nirenberg \cite{John-Nirenberg} established the following
deep property of functions from $BMO$.

\begin{theorem}
\label{Theorem}\cite{John-Nirenberg} If $b\in BMO({\mathbb{R}^{n}})$ and
$B\left(  x,r\right)  $ is a ball, then%
\[
\left \vert \left \{  x\in B\left(  x,r\right)  \,:\,|b(x)-b_{B\left(
x,r\right)  }|>\xi \right \}  \right \vert \leq|B\left(  x,r\right)  |\exp \left(
-\frac{\xi}{C\Vert b\Vert_{\ast}}\right)  ,~~~\xi>0,
\]
where $C$ depends only on the dimension $n$.
\end{theorem}

By Theorem \ref{Theorem}, we can get the following results.

\begin{corollary}
\cite{John-Nirenberg} Let $b\in BMO({\mathbb{R}^{n}})$. Then, for any $q>1$,%
\begin{equation}
\Vert b\Vert_{\ast}\thickapprox \sup_{x\in{\mathbb{R}^{n}},r>0}\left(  \frac
{1}{|B(x,r)|}%
{\displaystyle \int \limits_{B(x,r)}}
|b(y)-b_{B(x,r)}|^{p}dy\right)  ^{\frac{1}{p}}\label{5.1}%
\end{equation}
is valid.
\end{corollary}

\begin{corollary}
\label{Corollary} Let $b\in BMO({\mathbb{R}^{n}})$. Then there is a constant
$C>0$ such that
\begin{equation}
\left \vert b_{B(x,r)}-b_{B(x,t)}\right \vert \leq C\Vert b\Vert_{\ast}\left(
1+\ln \frac{t}{r}\right)  \text{ }~\text{for}~0<2r<t,\label{5.2}%
\end{equation}
and for any $q>1$, it is easy to see that%
\begin{equation}
\left \Vert b-\left(  b\right)  _{B}\right \Vert _{L_{q}\left(  B\right)  }\leq
Cr^{\frac{n}{q}}\Vert b\Vert_{\ast}\left(  1+\ln \frac{t}{r}\right)  .\label{c}%
\end{equation}
where $C$ is independent of $b$, $x$, $r$ and $t$.
\end{corollary}

\section{Proofs of the main results}

\subsection{Proof of Theorem \ref{teo9*}}

In the proof we have used the idea in \cite{Gao-Tang}. It suffices to show
that%
\[
\mu_{j,\Omega,b}^{L}f\left(  x\right)  \leq \mu_{j,\Omega,b}f\left(  x\right)
+CM_{\Omega,b}f\left(  x\right)  ,\text{ a.e. }x\in{\mathbb{R}^{n},}%
\]
where $M_{\Omega,b}$ denotes commutator of the Hardy-Littlewood maximal
operator with rough kernel.

Fix $x\in{\mathbb{R}^{n}}$ and let $r=\rho \left(  x\right)  $.%
\begin{align*}
\mu_{j,\Omega,b}^{L}f\left(  x\right)   & \leq \left(
{\displaystyle \int \limits_{0}^{r}}
\left \vert
{\displaystyle \int \limits_{\left \vert x-y\right \vert \leq t}}
\left \vert \Omega \left(  x-y\right)  \right \vert K_{j}^{L}\left(  x,y\right)
\left[  b\left(  x\right)  -b\left(  y\right)  \right]  f\left(  y\right)
dy\right \vert ^{2}\frac{dt}{t^{3}}\right)  ^{\frac{1}{2}}\\
& +\left(
{\displaystyle \int \limits_{r}^{\infty}}
\left \vert
{\displaystyle \int \limits_{\left \vert x-y\right \vert \leq r}}
\left \vert \Omega \left(  x-y\right)  \right \vert K_{j}^{L}\left(  x,y\right)
\left[  b\left(  x\right)  -b\left(  y\right)  \right]  f\left(  y\right)
dy\right \vert ^{2}\frac{dt}{t^{3}}\right)  ^{\frac{1}{2}}\\
& +\left(
{\displaystyle \int \limits_{r}^{\infty}}
\left \vert
{\displaystyle \int \limits_{r<\left \vert x-y\right \vert \leq t}}
\left \vert \Omega \left(  x-y\right)  \right \vert K_{j}^{L}\left(  x,y\right)
\left[  b\left(  x\right)  -b\left(  y\right)  \right]  f\left(  y\right)
dy\right \vert ^{2}\frac{dt}{t^{3}}\right)  ^{\frac{1}{2}}\\
& \leq \left(
{\displaystyle \int \limits_{0}^{r}}
\left \vert
{\displaystyle \int \limits_{\left \vert x-y\right \vert \leq t}}
\left \vert \Omega \left(  x-y\right)  \right \vert \left \vert K_{j}^{L}\left(
x,y\right)  -K_{j}\left(  x,y\right)  \right \vert \left[  b\left(  x\right)
-b\left(  y\right)  \right]  f\left(  y\right)  dy\right \vert ^{2}\frac
{dt}{t^{3}}\right)  ^{\frac{1}{2}}\\
& +\left(
{\displaystyle \int \limits_{0}^{r}}
\left \vert
{\displaystyle \int \limits_{\left \vert x-y\right \vert \leq t}}
\left \vert \Omega \left(  x-y\right)  \right \vert K_{j}\left(  x,y\right)
\left[  b\left(  x\right)  -b\left(  y\right)  \right]  f\left(  y\right)
dy\right \vert ^{2}\frac{dt}{t^{3}}\right)  ^{\frac{1}{2}}\\
& +\left(
{\displaystyle \int \limits_{r}^{\infty}}
\left \vert
{\displaystyle \int \limits_{\left \vert x-y\right \vert \leq r}}
\left \vert \Omega \left(  x-y\right)  \right \vert K_{j}^{L}\left(  x,y\right)
\left[  b\left(  x\right)  -b\left(  y\right)  \right]  f\left(  y\right)
dy\right \vert ^{2}\frac{dt}{t^{3}}\right)  ^{\frac{1}{2}}\\
& +\left(
{\displaystyle \int \limits_{r}^{\infty}}
\left \vert
{\displaystyle \int \limits_{r<\left \vert x-y\right \vert \leq t}}
\left \vert \Omega \left(  x-y\right)  \right \vert K_{j}^{L}\left(  x,y\right)
\left[  b\left(  x\right)  -b\left(  y\right)  \right]  f\left(  y\right)
dy\right \vert ^{2}\frac{dt}{t^{3}}\right)  ^{\frac{1}{2}}\\
& =E_{1}+E_{2}+E_{3}+E_{4}.
\end{align*}

For $E_{1}$, by Lemma \ref{lemma1}, we have%
\begin{align*}
E_{1}  & \leq C\left(
{\displaystyle \int \limits_{0}^{r}}
\left \vert
{\displaystyle \int \limits_{\left \vert x-y\right \vert \leq t}}
\frac{1}{\left \vert x-y\right \vert ^{n-1}}\left(  \frac{\left \vert
x-y\right \vert }{\rho \left(  x\right)  }\right)  ^{\delta}\left \vert
\Omega \left(  x-y\right)  \right \vert \left[  b\left(  x\right)  -b\left(
y\right)  \right]  f\left(  y\right)  dy\right \vert ^{2}\frac{dt}{t^{3}%
}\right)  ^{\frac{1}{2}}\\
& \leq Cr^{-\delta}\left(
{\displaystyle \int \limits_{0}^{r}}
\left \vert
{\displaystyle \sum \limits_{k=-\infty}^{0}}
\frac{1}{\left(  2^{k-1}t\right)  ^{n-\delta-1}}%
{\displaystyle \int \limits_{\left \vert x-y\right \vert \leq2^{k}t}}
\left \vert \Omega \left(  x-y\right)  \right \vert \left[  b\left(  x\right)
-b\left(  y\right)  \right]  f\left(  y\right)  dy\right \vert ^{2}\frac
{dt}{t^{3}}\right)  ^{\frac{1}{2}}\\
& \leq Cr^{-\delta}\left(
{\displaystyle \int \limits_{0}^{r}}
\left \vert
{\displaystyle \sum \limits_{k=-\infty}^{0}}
\frac{\left(  2^{k}\right)  ^{\delta+1}t^{\delta+1}}{\left(  2^{k}t\right)
^{n}}%
{\displaystyle \int \limits_{\left \vert x-y\right \vert \leq2^{k}t}}
\left \vert \Omega \left(  x-y\right)  \right \vert \left[  b\left(  x\right)
-b\left(  y\right)  \right]  f\left(  y\right)  dy\right \vert ^{2}\frac
{dt}{t^{3}}\right)  ^{\frac{1}{2}}\\
& \leq Cr^{-\delta}\left(
{\displaystyle \int \limits_{0}^{r}}
\left \vert
{\displaystyle \sum \limits_{k=-\infty}^{0}}
\left(  2^{k}\right)  ^{\delta+1}t^{\delta+1}M_{\Omega,b}f\left(  x\right)
\right \vert ^{2}\frac{dt}{t^{3}}\right)  ^{\frac{1}{2}}\\
& \leq Cr^{-\delta}\left(
{\displaystyle \int \limits_{0}^{r}}
t^{2\delta-1}dt\right)  ^{\frac{1}{2}}M_{\Omega,b}f\left(  x\right) \\
& \leq CM_{\Omega,b}f\left(  x\right)  .
\end{align*}

Obviously,%
\[
E_{2}\leq \mu_{j,\Omega,b}f\left(  x\right)  .
\]

For $E_{3}$, using Lemma \ref{lemma1} again, we get%
\begin{align*}
E_{3}  & \leq C\left(
{\displaystyle \int \limits_{r}^{\infty}}
\left \vert
{\displaystyle \int \limits_{\left \vert x-y\right \vert \leq r}}
\frac{1}{\left \vert x-y\right \vert ^{n-1}}\left \vert \Omega \left(  x-y\right)
\right \vert \left[  b\left(  x\right)  -b\left(  y\right)  \right]  f\left(
y\right)  dy\right \vert ^{2}\frac{dt}{t^{3}}\right)  ^{\frac{1}{2}}\\
& \leq C\left(
{\displaystyle \int \limits_{r}^{\infty}}
\left \vert
{\displaystyle \sum \limits_{k=-\infty}^{0}}
\frac{1}{\left(  2^{k-1}r\right)  ^{n-1}}%
{\displaystyle \int \limits_{\left \vert x-y\right \vert \leq2^{k}r}}
\left \vert \Omega \left(  x-y\right)  \right \vert \left[  b\left(  x\right)
-b\left(  y\right)  \right]  f\left(  y\right)  dy\right \vert ^{2}\frac
{dt}{t^{3}}\right)  ^{\frac{1}{2}}\\
& \leq C\left(
{\displaystyle \int \limits_{r}^{\infty}}
\left \vert
{\displaystyle \sum \limits_{k=-\infty}^{0}}
\frac{2^{k}r}{\left(  2^{k}r\right)  ^{n}}%
{\displaystyle \int \limits_{\left \vert x-y\right \vert \leq2^{k}r}}
\left \vert \Omega \left(  x-y\right)  \right \vert \left[  b\left(  x\right)
-b\left(  y\right)  \right]  f\left(  y\right)  dy\right \vert ^{2}\frac
{dt}{t^{3}}\right)  ^{\frac{1}{2}}\\
& \leq C\left(
{\displaystyle \int \limits_{r}^{\infty}}
\left \vert
{\displaystyle \sum \limits_{k=-\infty}^{0}}
2^{k}rM_{\Omega,b}f\left(  x\right)  \right \vert ^{2}\frac{dt}{t^{3}}\right)
^{\frac{1}{2}}\\
& \leq Cr\left(
{\displaystyle \int \limits_{r}^{\infty}}
\frac{dt}{t^{3}}\right)  ^{\frac{1}{2}}M_{\Omega,b}f\left(  x\right) \\
& \leq CM_{\Omega,b}f\left(  x\right)  .
\end{align*}

It remains to estimate $E_{4}$. By Lemma \ref{lemma1}, we obtain%
\begin{align*}
E_{4}  & \leq C\left(
{\displaystyle \int \limits_{r}^{\infty}}
\left \vert r%
{\displaystyle \int \limits_{r<\left \vert x-y\right \vert \leq t}}
\left \vert \Omega \left(  x-y\right)  \right \vert \left[  b\left(  x\right)
-b\left(  y\right)  \right]  \frac{\left \vert f\left(  y\right)  \right \vert
}{\left \vert x-y\right \vert ^{n}}dy\right \vert ^{2}\frac{dt}{t^{3}}\right)
^{\frac{1}{2}}\\
& \leq Cr\left(
{\displaystyle \int \limits_{r}^{\infty}}
\left \vert
{\displaystyle \sum \limits_{k=0}^{ \left[  \log_{2}t/r\right]  +1}}
\left(  2^{k}r\right)  ^{-n}%
{\displaystyle \int \limits_{\left \vert x-y\right \vert \leq2^{k}r}}
\left \vert \Omega \left(  x-y\right)  \right \vert \left[  b\left(  x\right)
-b\left(  y\right)  \right]  f\left(  y\right)  dy\right \vert ^{2}\frac
{dt}{t^{3}}\right)  ^{\frac{1}{2}}\\
& \leq Cr\left(
{\displaystyle \int \limits_{r}^{\infty}}
\left \vert \left(  \left[  \log_{2}\frac{t}{r}\right]  +1\right)  M_{\Omega
,b}f\left(  x\right)  \right \vert ^{2}\frac{dt}{t^{3}}\right)  ^{\frac{1}{2}%
}\\
& \leq Cr\left(
{\displaystyle \int \limits_{r}^{\infty}}
\frac{t}{r}M_{\Omega,b}f\left(  x\right)  ^{2}\frac{dt}{t^{3}}\right)
^{\frac{1}{2}}\\
& \leq CM_{\Omega,b}f\left(  x\right)  .
\end{align*}

Thus, Theorem \ref{teo9*} is proved.

\subsection{Proof of Theorem \ref{teo10*}}

For $x\in B\left(  x_{0},t\right)  $, notice that $\Omega$ is homogenous of
degree zero and $\Omega \in L_{q}(S^{n-1})$, $1<q\leq \infty$. Then, we obtain%
\begin{align}
\left(  \int \limits_{B\left(  x_{0},t\right)  }\left \vert \Omega \left(
x-y\right)  \right \vert ^{q}dy\right)  ^{\frac{1}{q}}  & =\left(
\int \limits_{B\left(  x-x_{0},t\right)  }\left \vert \Omega \left(  z\right)
\right \vert ^{q}dz\right)  ^{\frac{1}{q}}\nonumber \\
& \leq \left(  \int \limits_{B\left(  0,t+\left \vert x-x_{0}\right \vert \right)
}\left \vert \Omega \left(  z\right)  \right \vert ^{q}dz\right)  ^{\frac{1}{q}%
}\nonumber \\
& \leq \left(  \int \limits_{B\left(  0,2t\right)  }\left \vert \Omega \left(
z\right)  \right \vert ^{q}dz\right)  ^{\frac{1}{q}}\nonumber \\
& =\left(  \int \limits_{0}^{2t}\int \limits_{S^{n-1}}\left \vert \Omega \left(
z^{\prime}\right)  \right \vert ^{q}d\sigma \left(  z^{\prime}\right)
r^{n-1}dr\right)  ^{\frac{1}{q}}\nonumber \\
& =C\left \Vert \Omega \right \Vert _{L_{q}\left(  S^{n-1}\right)  }\left \vert
B\left(  x_{0},2t\right)  \right \vert ^{\frac{1}{q}}.\label{10}%
\end{align}

Let $1<p<\infty$ and $q^{\prime}\leq p$. For any $x_{0}\in{\mathbb{R}^{n}} $,
set $B=B\left(  x_{0},r\right)  $ for the ball centered at $x_{0}$ and of
radius $r$ and $2B=B\left(  x_{0},2r\right)  $. We represent $f$ as%
\[
f=f_{1}+f_{2},\text{ \  \ }f_{1}\left(  y\right)  =f\left(  y\right)  \chi
_{2B}\left(  y\right)  ,\text{ \  \ }f_{2}\left(  y\right)  =f\left(  y\right)
\chi_{\, \! \left(  2B\right)  ^{C}}\left(  y\right)  ,\text{ \  \ }r>0
\]
and have%
\[
\left \Vert \mu_{j,\Omega,b}^{L}f\right \Vert _{L_{p}\left(  B\right)  }%
\leq \left \Vert \mu_{j,\Omega,b}^{L}f_{1}\right \Vert _{L_{p}\left(  B\right)
}+\left \Vert \mu_{j,\Omega,b}^{L}f_{2}\right \Vert _{L_{p}\left(  B\right)  }.
\]

Since $f_{1}\in$ $L_{p}({\mathbb{R}^{n}})$, $\mu_{j,\Omega,b}^{L}f_{1}\in$
$L_{p}({\mathbb{R}^{n}})$, from the boundedness of $\mu_{j,\Omega,b}^{L}$ on
$L_{p}({\mathbb{R}^{n}})$ (see Theorem \ref{teo9*}) it follows that:%
\begin{align*}
\left \Vert \mu_{j,\Omega,b}^{L}f_{1}\right \Vert _{L_{p}\left(  B\right)  }  &
\leq \left \Vert \mu_{j,\Omega,b}^{L}f_{1}\right \Vert _{L_{p}\left(
{\mathbb{R}^{n}}\right)  }\\
& \lesssim \left \Vert b\right \Vert _{\ast}\left \Vert f_{1}\right \Vert
_{L_{p}\left(  {\mathbb{R}^{n}}\right)  }=\left \Vert b\right \Vert _{\ast
}\left \Vert f\right \Vert _{L_{p}\left(  2B\right)  }.
\end{align*}

It is known that $x\in B$, $y\in \left(  2B\right)  ^{C}$, which implies
$\frac{1}{2}\left \vert x_{0}-y\right \vert \leq \left \vert x-y\right \vert
\leq \frac{3}{2}\left \vert x_{0}-y\right \vert $. Then for $x\in B$, we have%
\begin{align*}
\left \vert \mu_{j,\Omega,b}^{L}f_{2}\left(  x\right)  \right \vert  &
\lesssim \int \limits_{{\mathbb{R}^{n}}}\frac{\left \vert \Omega \left(
x-y\right)  \right \vert }{\left \vert x-y\right \vert ^{n}}\left \vert b\left(
y\right)  -b\left(  x\right)  \right \vert \left \vert f\left(  y\right)
\right \vert dy\\
& \approx \int \limits_{\left(  2B\right)  ^{C}}\frac{\left \vert \Omega \left(
x-y\right)  \right \vert }{\left \vert x_{0}-y\right \vert ^{n}}\left \vert
b\left(  y\right)  -b\left(  x\right)  \right \vert \left \vert f\left(
y\right)  \right \vert dy.
\end{align*}

Hence we get%
\begin{align*}
\left \Vert \mu_{j,\Omega,b}^{L}f_{2}\right \Vert _{L_{p}\left(  B\right)  }  &
\lesssim \left(  \int \limits_{B}\left(  \int \limits_{\left(  2B\right)  ^{C}%
}\frac{\left \vert \Omega \left(  x-y\right)  \right \vert }{\left \vert
x_{0}-y\right \vert ^{n}}\left \vert b\left(  y\right)  -b\left(  x\right)
\right \vert \left \vert f\left(  y\right)  \right \vert dy\right)
^{p}dx\right)  ^{\frac{1}{p}}\\
& \lesssim \left(  \int \limits_{B}\left(  \int \limits_{\left(  2B\right)  ^{C}%
}\frac{\left \vert \Omega \left(  x-y\right)  \right \vert }{\left \vert
x_{0}-y\right \vert ^{n}}\left \vert b\left(  y\right)  -b_{B}\right \vert
\left \vert f\left(  y\right)  \right \vert dy\right)  ^{p}dx\right)  ^{\frac
{1}{p}}\\
& +\left(  \int \limits_{B}\left(  \int \limits_{\left(  2B\right)  ^{C}}%
\frac{\left \vert \Omega \left(  x-y\right)  \right \vert }{\left \vert
x_{0}-y\right \vert ^{n}}\left \vert b\left(  x\right)  -b_{B}\right \vert
\left \vert f\left(  y\right)  \right \vert dy\right)  ^{p}dx\right)  ^{\frac
{1}{p}}\\
& =J_{1}+J_{2}.
\end{align*}

We have the following estimation of $J_{1}$. When $q^{\prime}\leq p$ and
$\frac{1}{\mu}+\frac{1}{p}+\frac{1}{q}=1$, by the Fubini's theorem%
\begin{align*}
J_{1}  & \approx r^{\frac{n}{p}}\int \limits_{\left(  2B\right)  ^{C}}%
\frac{\left \vert \Omega \left(  x-y\right)  \right \vert }{\left \vert
x_{0}-y\right \vert ^{n}}\left \vert b\left(  y\right)  -b_{B}\right \vert
\left \vert f\left(  y\right)  \right \vert dy\\
& \approx r^{\frac{n}{p}}\int \limits_{\left(  2B\right)  ^{C}}\left \vert
\Omega \left(  x-y\right)  \right \vert \left \vert b\left(  y\right)
-b_{B}\right \vert \left \vert f\left(  y\right)  \right \vert \int
\limits_{\left \vert x_{0}-y\right \vert }^{\infty}\frac{dt}{t^{n+1}}dy\\
& \approx r^{\frac{n}{p}}\int \limits_{2r}^{\infty}\int \limits_{2r\leq
\left \vert x_{0}-y\right \vert \leq t}\left \vert \Omega \left(  x-y\right)
\right \vert \left \vert b\left(  y\right)  -b_{B}\right \vert \left \vert
f\left(  y\right)  \right \vert dy\frac{dt}{t^{n+1}}\\
& \lesssim r^{\frac{n}{p}}\int \limits_{2r}^{\infty}\int \limits_{B\left(
x_{0},t\right)  }\left \vert \Omega \left(  x-y\right)  \right \vert \left \vert
b\left(  y\right)  -b_{B}\right \vert \left \vert f\left(  y\right)  \right \vert
dy\frac{dt}{t^{n+1}}\text{ holds.}%
\end{align*}

Applying the H\"{o}lder's inequality and by (\ref{10}), (\ref{5.1}),
(\ref{5.2}) and (\ref{c}), we get%
\begin{align*}
J_{1}  & \lesssim r^{\frac{n}{p}}\int \limits_{2r}^{\infty}\int
\limits_{B\left(  x_{0},t\right)  }\left \vert \Omega \left(  x-y\right)
\right \vert \left \vert b\left(  y\right)  -b_{B\left(  x_{0},t\right)
}\right \vert \left \vert f\left(  y\right)  \right \vert dy\frac{dt}{t^{n+1}}\\
& +r^{\frac{n}{p}}\int \limits_{2r}^{\infty}\left \vert b_{B\left(
x_{0},r\right)  }-b_{B\left(  x_{0},t\right)  }\right \vert \int
\limits_{B\left(  x_{0},t\right)  }\left \vert \Omega \left(  x-y\right)
\right \vert \left \vert f\left(  y\right)  \right \vert dy\frac{dt}{t^{n+1}}\\
& \lesssim r^{\frac{n}{p}}\int \limits_{2r}^{\infty}\left \Vert \Omega \left(
\cdot-y\right)  \right \Vert _{L_{q}\left(  B\left(  x_{0},t\right)  \right)
}\left \Vert \left(  b\left(  \cdot \right)  -b_{B\left(  x_{0},t\right)
}\right)  \right \Vert _{L_{\mu}\left(  B\left(  x_{0},t\right)  \right)
}\left \Vert f\right \Vert _{L_{p}\left(  B\left(  x_{0},t\right)  \right)
}\frac{dt}{t^{n+1}}\\
& +r^{\frac{n}{p}}\int \limits_{2r}^{\infty}\left \vert b_{B\left(
x_{0},r\right)  }-b_{B\left(  x_{0},t\right)  }\right \vert \left \Vert
\Omega \left(  \cdot-y\right)  \right \Vert _{L_{q}\left(  B\left(
x_{0},t\right)  \right)  }\left \Vert f\right \Vert _{L_{p}\left(  B\left(
x_{0},t\right)  \right)  }\left \vert B\left(  x_{0},t\right)  \right \vert
^{1-\frac{1}{p}-\frac{1}{q}}\frac{dt}{t^{n+1}}\\
& \lesssim \Vert b\Vert_{\ast}\,r^{\frac{n}{p}}\int \limits_{2r}^{\infty}\left(
1+\ln \frac{t}{r}\right)  \left \Vert f\right \Vert _{L_{p}\left(  B\left(
x_{0},t\right)  \right)  }\frac{dt}{t^{\frac{n}{p}+1}}.
\end{align*}

In order to estimate $J_{2}$ note that%
\[
J_{2}=\left \Vert \left(  b\left(  \cdot \right)  -b_{B\left(  x_{0},t\right)
}\right)  \right \Vert _{L_{p}\left(  B\left(  x_{0},t\right)  \right)  }%
\int \limits_{\left(  2B\right)  ^{C}}\frac{\left \vert \Omega \left(
x-y\right)  \right \vert }{\left \vert x_{0}-y\right \vert ^{n}}\left \vert
f\left(  y\right)  \right \vert dy.
\]

By (\ref{5.1}), we get%
\[
J_{2}\lesssim \Vert b\Vert_{\ast}\,r^{\frac{n}{p}}\int \limits_{\left(
2B\right)  ^{C}}\frac{\left \vert \Omega \left(  x-y\right)  \right \vert
}{\left \vert x_{0}-y\right \vert ^{n}}\left \vert f\left(  y\right)  \right \vert
dy.
\]

Applying the H\"{o}lder's inequality, we get%
\[
J_{2}\lesssim \Vert b\Vert_{\ast}\,r^{\frac{n}{p}}\int \limits_{2r}^{\infty
}\left \Vert f\right \Vert _{L_{p}\left(  B\left(  x_{0},t\right)  \right)
}\left \Vert \Omega \left(  x-\cdot \right)  \right \Vert _{L_{q}\left(  B\left(
x_{0},t\right)  \right)  }\left \vert B\left(  x_{0},t\right)  \right \vert
^{1-\frac{1}{p}-\frac{1}{q}}\frac{dt}{t^{n+1}}.
\]

Thus, by (\ref{10}) we get%
\[
J_{2}\lesssim \Vert b\Vert_{\ast}\,r^{\frac{n}{p}}\int \limits_{2r}^{\infty
}\left \Vert f\right \Vert _{L_{p}\left(  B\left(  x_{0},t\right)  \right)
}\frac{dt}{t^{\frac{n}{p}+1}}.
\]

Summing up $J_{1}$ and $J_{2}$, for all $p\in \left(  1,\infty \right)  $ we get%
\[
\left \Vert \mu_{j,\Omega,b}^{L}f_{2}\right \Vert _{L_{p}\left(  B\right)
}\lesssim \Vert b\Vert_{\ast}\,r^{\frac{n}{p}}\int \limits_{2r}^{\infty}\left(
1+\ln \frac{t}{r}\right)  \left \Vert f\right \Vert _{L_{p}\left(  B\left(
x_{0},t\right)  \right)  }\frac{dt}{t^{\frac{n}{p}+1}}.
\]

Finally, we have the following%
\[
\left \Vert \mu_{j,\Omega,b}^{L}f\right \Vert _{L_{p}\left(  B\right)  }%
\lesssim \left \Vert b\right \Vert _{\ast}\left \Vert f\right \Vert _{L_{p}\left(
2B\right)  }+\Vert b\Vert_{\ast}\,r^{\frac{n}{p}}\int \limits_{2r}^{\infty
}\left(  1+\ln \frac{t}{r}\right)  \left \Vert f\right \Vert _{L_{p}\left(
B\left(  x_{0},t\right)  \right)  }\frac{dt}{t^{\frac{n}{p}+1}}.
\]
On the other hand, we have%

\begin{align}
\left \Vert f\right \Vert _{L_{p}\left(  2B\right)  }  & \approx r^{\frac{n}{p}%
}\left \Vert f\right \Vert _{L_{p}\left(  2B\right)  }\int \limits_{2r}^{\infty
}\frac{dt}{t^{\frac{n}{p}+1}}\nonumber \\
& \leq r^{\frac{n}{p}}\int \limits_{2r}^{\infty}\left \Vert f\right \Vert
_{L_{p}\left(  B\left(  x_{0},t\right)  \right)  }\frac{dt}{t^{\frac{n}{p}+1}%
}.\label{e313}%
\end{align}

By combining the above inequalities, we obtain%
\[
\Vert \mu_{j,\Omega,b}^{L}f\Vert_{L_{p}(B(x_{0},r))}\lesssim \Vert b\Vert_{\ast
}\,r^{\frac{n}{p}}\int \limits_{2r}^{\infty}\left(  1+\ln \frac{t}{r}\right)
t^{-\frac{n}{p}-1}\Vert f\Vert_{L_{p}(B(x_{0},t))}dt,
\]
which completes the proof of first statement.

Similarly to (\ref{10}), when $y\in B\left(  x_{0},t\right)  $, it is true that%

\begin{equation}
\left(  \int \limits_{B\left(  x_{0},r\right)  }\left \vert \Omega \left(
x-y\right)  \right \vert ^{q}dy\right)  ^{\frac{1}{q}}\leq C\left \Vert
\Omega \right \Vert _{L_{q}\left(  S^{n-1}\right)  }\left \vert B\left(
x_{0},\frac{3}{2}t\right)  \right \vert ^{\frac{1}{q}}.\label{314}%
\end{equation}

On the other hand when $p<q$, by the Fubini's theorem and the Minkowski
inequality, we get%
\begin{align*}
J_{1}  & \lesssim \left(  \int \limits_{B}\left \vert \int \limits_{2r}^{\infty
}\int \limits_{B\left(  x_{0},t\right)  }\left \vert b\left(  y\right)
-b_{B\left(  x_{0},t\right)  }\right \vert \left \vert f\left(  y\right)
\right \vert \left \vert \Omega \left(  x-y\right)  \right \vert dy\frac
{dt}{t^{n+1}}\right \vert ^{p}dx\right)  ^{\frac{1}{p}}\\
& +\left(  \int \limits_{B}\left \vert \int \limits_{2r}^{\infty}\left \vert
b_{B\left(  x_{0},r\right)  }-b_{B\left(  x_{0},t\right)  }\right \vert
\int \limits_{B\left(  x_{0},t\right)  }\left \vert f\left(  y\right)
\right \vert \left \vert \Omega \left(  x-y\right)  \right \vert dy\frac
{dt}{t^{n+1}}\right \vert ^{p}dx\right)  ^{\frac{1}{p}}\\
& \lesssim \int \limits_{2r}^{\infty}\int \limits_{B\left(  x_{0},t\right)
}\left \vert b\left(  y\right)  -b_{B\left(  x_{0},t\right)  }\right \vert
\left \vert f\left(  y\right)  \right \vert \left \Vert \Omega \left(
\cdot-y\right)  \right \Vert _{L_{p}\left(  B\left(  x_{0},t\right)  \right)
}dy\frac{dt}{t^{n+1}}\\
& +\int \limits_{2r}^{\infty}\left \vert b_{B\left(  x_{0},r\right)
}-b_{B\left(  x_{0},t\right)  }\right \vert \int \limits_{B\left(
x_{0},t\right)  }\left \vert f\left(  y\right)  \right \vert \left \Vert
\Omega \left(  \cdot-y\right)  \right \Vert _{L_{p}\left(  B\left(
x_{0},t\right)  \right)  }dy\frac{dt}{t^{n+1}}\\
& \lesssim \left \vert B\right \vert ^{\frac{1}{p}-\frac{1}{q}}\int
\limits_{2r}^{\infty}\int \limits_{B\left(  x_{0},t\right)  }\left \vert
b\left(  y\right)  -b_{B\left(  x_{0},t\right)  }\right \vert \left \vert
f\left(  y\right)  \right \vert \left \Vert \Omega \left(  \cdot-y\right)
\right \Vert _{L_{q}\left(  B\left(  x_{0},t\right)  \right)  }dy\frac
{dt}{t^{n+1}}\\
& +\left \vert B\right \vert ^{\frac{1}{p}-\frac{1}{q}}\int \limits_{2r}^{\infty
}\left \vert b_{B\left(  x_{0},r\right)  }-b_{B\left(  x_{0},t\right)
}\right \vert \int \limits_{B\left(  x_{0},t\right)  }\left \vert f\left(
y\right)  \right \vert \left \Vert \Omega \left(  \cdot-y\right)  \right \Vert
_{L_{q}\left(  B\left(  x_{0},t\right)  \right)  }dy\frac{dt}{t^{n+1}}.
\end{align*}

Applying the H\"{o}lder's inequality and by (\ref{314}), (\ref{5.1}),
(\ref{5.2}) and Lemma \ref{lemma 100}, we get%
\begin{align*}
J_{1}  & \lesssim r^{\frac{n}{p}-\frac{n}{q}}\int \limits_{2r}^{\infty
}\left \Vert \left(  b\left(  \cdot \right)  -b_{B\left(  x_{0},t\right)
}\right)  f\right \Vert _{L_{1}\left(  B\left(  x_{0},t\right)  \right)
}\left \vert B\left(  x_{0},\frac{3}{2}t\right)  \right \vert ^{\frac{1}{q}%
}\frac{dt}{t^{n+1}}\\
& +r^{\frac{n}{p}-\frac{n}{q}}\int \limits_{2r}^{\infty}\left \vert b_{B\left(
x_{0},r\right)  }-b_{B\left(  x_{0},t\right)  }\right \vert \left \Vert
f\right \Vert _{L_{p}\left(  B\left(  x_{0},t\right)  \right)  }\left \vert
B\left(  x_{0},\frac{3}{2}t\right)  \right \vert ^{\frac{1}{q}}\frac
{dt}{t^{\frac{n}{p}+1}}\\
& \lesssim r^{\frac{n}{p}-\frac{n}{q}}\int \limits_{2r}^{\infty}\left \Vert
\left(  b\left(  \cdot \right)  -b_{B\left(  x_{0},t\right)  }\right)
\right \Vert _{L_{p^{\prime}}\left(  B\left(  x_{0},t\right)  \right)
}\left \Vert f\right \Vert _{L_{p}\left(  B\left(  x_{0},t\right)  \right)
}t^{\frac{n}{q}}\frac{dt}{t^{n+1}}\\
& +r^{\frac{n}{p}-\frac{n}{q}}\int \limits_{2r}^{\infty}\left \vert b_{B\left(
x_{0},r\right)  }-b_{B\left(  x_{0},t\right)  }\right \vert \left \Vert
f\right \Vert _{L_{p}\left(  B\left(  x_{0},t\right)  \right)  }t^{\frac{n}{q}%
}\frac{dt}{t^{\frac{n}{p}+1}}\\
& \lesssim \Vert b\Vert_{\ast}\,r^{\frac{n}{p}-\frac{n}{q}}\int \limits_{2r}%
^{\infty}\left(  1+\ln \frac{t}{r}\right)  t^{\frac{n}{q}-\frac{n}{p}%
-1}\left \Vert f\right \Vert _{L_{p}\left(  B\left(  x_{0},t\right)  \right)
}dt.
\end{align*}

Let $\frac{1}{p}=\frac{1}{\nu}+\frac{1}{q}$, then for $J_{2}$, by the Fubini's
theorem, the Minkowski inequality, the H\"{o}lder's inequality and from
(\ref{314}), we get%
\begin{align*}
J_{2}  & \lesssim \left(  \int \limits_{B}\left \vert \int \limits_{2r}^{\infty
}\int \limits_{B\left(  x_{0},t\right)  }\left \vert f\left(  y\right)
\right \vert \left \vert b\left(  x\right)  -b_{B}\right \vert \left \vert
\Omega \left(  x-y\right)  \right \vert dy\frac{dt}{t^{n+1}}\right \vert
^{p}dx\right)  ^{\frac{1}{p}}\\
& \lesssim \int \limits_{2r}^{\infty}\int \limits_{B\left(  x_{0},t\right)
}\left \vert f\left(  y\right)  \right \vert \left \Vert \left(  b\left(
\cdot \right)  -b_{B}\right)  \Omega \left(  \cdot-y\right)  \right \Vert
_{L_{p}\left(  B\right)  }dy\frac{dt}{t^{n+1}}\\
& \lesssim \int \limits_{2r}^{\infty}\int \limits_{B\left(  x_{0},t\right)
}\left \vert f\left(  y\right)  \right \vert \left \Vert b\left(  \cdot \right)
-b_{B}\right \Vert _{L_{\nu}\left(  B\right)  }\left \Vert \Omega \left(
\cdot-y\right)  \right \Vert _{L_{q}\left(  B\right)  }dy\frac{dt}{t^{n+1}}\\
& \lesssim \Vert b\Vert_{\ast}\left \vert B\right \vert ^{\frac{1}{p}-\frac{1}%
{q}}\int \limits_{2r}^{\infty}\int \limits_{B\left(  x_{0},t\right)  }\left \vert
f\left(  y\right)  \right \vert \left \Vert \Omega \left(  \cdot-y\right)
\right \Vert _{L_{q}\left(  B\right)  }dy\frac{dt}{t^{n+1}}\\
& \lesssim \Vert b\Vert_{\ast}r^{\frac{n}{p}-\frac{n}{q}}\int \limits_{2r}%
^{\infty}\left \Vert f\right \Vert _{L_{1}\left(  B\left(  x_{0},t\right)
\right)  }\left \vert B\left(  x_{0},\frac{3}{2}t\right)  \right \vert
^{\frac{1}{q}}\frac{dt}{t^{n+1}}\\
& \lesssim \Vert b\Vert_{\ast}\,r^{\frac{n}{p}-\frac{n}{q}}\int \limits_{2r}%
^{\infty}\left(  1+\ln \frac{t}{r}\right)  t^{\frac{n}{q}-\frac{n}{p}-1}\Vert
f\Vert_{L_{p}(B(x_{0},t))}dt.
\end{align*}

By combining the above estimates, we complete the proof of Theorem
\ref{teo10*}.

\subsection{Proof of Theorem \ref{teo11}}

The statement is derived from inequalities (\ref{47}) and (\ref{48}). Let
$q^{\prime}\leq p$. The estimation of the norm of the operator, that is, the
boundedness in the non-vanishing space follows from Theorem \ref{teo10*} and
condition (\ref{7*})
\begin{align*}
\Vert \mu_{j,\Omega,b}^{L}f\Vert_{VM_{p,\varphi_{2}}}  & =\sup \limits_{x\in
{\mathbb{R}^{n}},r>0}\varphi_{2}(x,r)^{-1}\Vert \mu_{j,\Omega,b}^{L}%
f\Vert_{L_{p}(B(x,r))}\\
& \lesssim \left \Vert b\right \Vert _{\ast}\sup_{x\in{\mathbb{R}^{n},}%
r>0}\varphi_{2}\left(  x,r\right)  ^{-1}r^{\frac{n}{p}}\int \limits_{r}%
^{\infty}\left(  1+\ln \frac{t}{r}\right)  \left \Vert f\right \Vert
_{L_{p}\left(  B\left(  x,t\right)  \right)  }\frac{dt}{t^{\frac{n}{p}+1}}\\
& \lesssim \left \Vert b\right \Vert _{\ast}\sup_{x\in{\mathbb{R}^{n},}%
r>0}\varphi_{2}\left(  x,r\right)  ^{-1}r^{\frac{n}{p}}\int \limits_{r}%
^{\infty}\left(  1+\ln \frac{t}{r}\right)  \varphi_{1}\left(  x,t\right)
\left[  \varphi_{1}\left(  x,t\right)  ^{-1}\left \Vert f\right \Vert
_{L_{p}\left(  B\left(  x,t\right)  \right)  }\right]  \frac{dt}{t^{\frac
{n}{p}+1}}\\
& \lesssim \left \Vert b\right \Vert _{\ast}\left \Vert f\right \Vert
_{VM_{p,\varphi_{1}}}\sup_{x\in{\mathbb{R}^{n},}r>0}\varphi_{2}\left(
x,r\right)  ^{-1}r^{\frac{n}{p}}\int \limits_{r}^{\infty}\left(  1+\ln \frac
{t}{r}\right)  \varphi_{1}\left(  x,t\right)  \frac{dt}{t^{\frac{n}{p}+1}}\\
& \lesssim \left \Vert b\right \Vert _{\ast}\left \Vert f\right \Vert
_{VM_{p,\varphi_{1}}}.
\end{align*}

So we only have to prove that%
\begin{equation}
\lim \limits_{r\rightarrow0}\sup \limits_{x\in{\mathbb{R}^{n}}}\varphi
_{1}(x,r)^{-1}\, \left \Vert f\right \Vert _{L_{p}\left(  B\left(  x,r\right)
\right)  }=0\text{ implies }\lim \limits_{r\rightarrow0}\sup \limits_{x\in
{\mathbb{R}^{n}}}\varphi_{2}(x,r)^{-1}\, \left \Vert \mu_{j,\Omega,b}%
^{L}f\right \Vert _{L_{p}\left(  B\left(  x,r\right)  \right)  }=0.\label{12*}%
\end{equation}

To show that $\sup \limits_{x\in{\mathbb{R}^{n}}}\varphi_{2}(x,r)^{-1}\,
\left \Vert \mu_{j,\Omega,b}^{L}f\right \Vert _{L_{p}\left(  B\left(
x,r\right)  \right)  }<\epsilon$ for small $r$, we split the right-hand side
of (\ref{47}):%
\begin{equation}
\varphi_{2}(x,r)^{-1}\, \left \Vert \mu_{j,\Omega,b}^{L}f\right \Vert
_{L_{p}\left(  B\left(  x,r\right)  \right)  }\leq C\left[  I_{\delta_{0}%
}\left(  x,r\right)  +J_{\delta_{0}}\left(  x,r\right)  \right]  ,\label{13*}%
\end{equation}
where $\delta_{0}>0$ (we may take $\delta_{0}<1$), and
\[
I_{\delta_{0}}\left(  x,r\right)  :=\left \Vert b\right \Vert _{\ast}%
\frac{r^{\frac{n}{p}}}{\varphi_{2}(x,r)}\left(
{\displaystyle \int \limits_{r}^{\delta_{0}}}
\left(  1+\ln \frac{t}{r}\right)  \varphi_{1}\left(  x,t\right)  t^{-\frac
{n}{p}-1}\left(  \varphi_{1}\left(  x,t\right)  ^{-1}\left \Vert f\right \Vert
_{L_{p}\left(  B\left(  x,t\right)  \right)  }\right)  dt\right)  ,
\]
and%
\[
J_{\delta_{0}}\left(  x,r\right)  :=\left \Vert b\right \Vert _{\ast}%
\frac{r^{\frac{n}{p}}}{\varphi_{2}(x,r)}\left(
{\displaystyle \int \limits_{\delta_{0}}^{\infty}}
\left(  1+\ln \frac{t}{r}\right)  \varphi_{1}\left(  x,t\right)  t^{-\frac
{n}{p}-1}\left(  \varphi_{1}\left(  x,t\right)  ^{-1}\left \Vert f\right \Vert
_{L_{p}\left(  B\left(  x,t\right)  \right)  }\right)  dt\right)
\]
and $r<\delta_{0}$. Now we choose any fixed $\delta_{0}>0$ such that%
\[
\sup \limits_{x\in{\mathbb{R}^{n}}}\varphi_{1}\left(  x,t\right)
^{-1}\left \Vert f\right \Vert _{L_{p}\left(  B\left(  x,t\right)  \right)
}<\frac{\epsilon}{2CC_{0}},
\]
where $C$ and $C_{0}$ are constants from (\ref{7*}) and (\ref{13*}), which is
possible since $f\in VM_{p,\varphi_{1}}$. This allows to estimate the first
term uniformly in $r\in \left(  0,\delta_{0}\right)  $:%
\[
\left \Vert b\right \Vert _{\ast}\sup \limits_{x\in{\mathbb{R}^{n}}}%
CI_{\delta_{0}}\left(  x,r\right)  <\frac{\epsilon}{2},\qquad0<r<\delta_{0}.
\]

The estimation of the second term may be obtained by choosing $r$ sufficiently
small. Indeed, by (\ref{2}) we have%
\[
J_{\delta_{0}}\left(  x,r\right)  \leq \left \Vert b\right \Vert _{\ast}%
c_{\delta_{0}}\left \Vert f\right \Vert _{VM_{p,\varphi}}\frac{r^{\frac{n}{p}}%
}{\varphi \left(  x,r\right)  },
\]
where $c_{\delta_{0}}$ is the constant from (\ref{6*}). Then, by (\ref{2}) it
suffices to choose $r$ small enough such that
\[
\sup \limits_{x\in{\mathbb{R}^{n}}}\frac{r^{\frac{n}{p}}}{\varphi(x,r)}%
\leq \frac{\epsilon}{2\left \Vert b\right \Vert _{\ast}c_{\delta_{0}}\left \Vert
f\right \Vert _{VM_{p,\varphi}}},
\]
which completes the proof of (\ref{12*}).

For the case of $p<q$, we can also use the same method, so we omit the
details. Thus, we obtain (\ref{10*}), which completes the proof of Theorem
\ref{teo11}.

\begin{remark}
Conditions (\ref{6*}) and (\ref{8*}) are not needed in the case when
$\varphi(x,r)$ does not depend on $x$, since (\ref{6*}) follows from
(\ref{7*}) and similarly, (\ref{8*}) follows from (\ref{9*}) in this case.
\end{remark}

\end{document}